\documentclass{elsart}

\usepackage{amssymb}
\usepackage{amsmath}

\newtheorem{trans}{Transformation}

\numberwithin{equation}{section}
\numberwithin{thm}{section}
 
\newcommand{\hgs}[6]{ {}_{#1}\phi_{#2} 
    \left[ \genfrac{}{}{0pt}{}{#3}{#4} ; {#5},{#6} \right]}
\newcommand{\W}[3]{ {}_{#1}W_{#2}\left( #3 \right) }
\newcommand{\binomial}[2]{ \genfrac{(}{)}{0pt}{}{#1}{#2} }

\begin{document}

\begin{frontmatter}
\journal{J. Math. Anal. Appl.}

\title{Identities of the \\Rogers-Ramanujan-Bailey Type}
\author{Andrew V. Sills}  
\address{Department of Mathematics, 
Rutgers University,
Hill Center--Busch Campus, Piscataway, NJ 08854; 
telephone 732-445-3488; fax 732-445-5530}
\ead{asills@math.rutgers.edu}
\ead[url]{http://www.math.rutgers.edu/\~{}asills}

\begin{abstract}
 A multiparameter generalization of the Bailey pair is defined in such a way as to
include as special cases all Bailey pairs considered by W. N. Bailey in his paper, 
``Identities of
the Rogers-Ramanujan type,"
\emph{Proc. London Math. Soc. (2)}, 50 (1949), 421--435.
This leads to the derivation of a number of elegant new Rogers-Ramanujan type
identities. 
\end{abstract}

\end{frontmatter}

\section{Introduction}
\subsection{Overview} \label{ov}
Recall the famous Rogers-Ramanujan identities:
\begin{thm}[The Rogers-Ramanujan Identities] 
 \begin{equation}\label{RRa1}
   \sum_{n=0}^\infty \frac{q^{n^2}}{(q;q)_n} = 
   \frac{(q^2, q^3, q^5; q^5)_\infty}{(q;q)_\infty},
 \end{equation} and
\begin{equation}\label{RRa2}
   \sum_{n=0}^\infty \frac{q^{n^2+n}}{(q;q)_n} = 
   \frac{(q, q^4, q^5; q^5)_\infty}{(q;q)_\infty},
 \end{equation}
where \[ (a;q)_m = \prod_{j=0}^{m-1} (1-aq^j), \]
      \[ (a;q)_\infty = \prod_{j=0}^\infty (1-aq^j), \] and
      \[ (a_1, a_2, \dots, a_r; q)_s = (a_1;q)_s (a_2;q)_s \dots (a_r;q)_s, \]
\end{thm} 
(Although the results in this paper may be considered purely from the
point of view of formal power series, they also yield identities of
analytic functions provided $|q|<1$.)

The Rogers-Ramanujan identities are due to L.~J.~Rogers~\cite{Rog1894},
and were rediscovered independently by S. Ramanujan~\cite{MacMahon} and
I. Schur~\cite{Schur}.  
Rogers (and later others) discovered many series--product 
identities similar in form to the Rogers-Ramanujan identities, and 
as such are referred to as ``identities of the Rogers-Ramanujan type."  
A number of Rogers-Ramanujan type identities were recorded by Ramanujan
in his Lost Notebook~\cite[Chapter 11]{rln}.
During World War II, W.~N.~Bailey undertook a thorough
study of Rogers' work connected with Rogers-Ramanujan type identities, 
and through the understanding he gained, was able to simplify and 
generalize Rogers' ideas in a pair of papers~(\cite{wnb1} and 
\cite{wnb2}). In the process, Bailey and Freeman Dyson (who served as
referee for Bailey's two papers~\cite[p. 14]{fjd}) discovered a 
number of new Rogers-Ramanujan type identities.  

   Bailey and his student L.J. Slater~(\cite{ljs1}, 
\cite{ljs2}) only
considered identities of single-fold series and infinite products.  
Bailey comments in passing
\cite[p. 4, \S 4]{wnb2} that ``the most general formulae for basic 
series (apart 
from those already given) are too involved to be of any great interest" and 
as such,
rejected multisums from consideration.  However, G.E. Andrews reversed this 
prejudice
against multisum Rogers-Ramanujan type identities by presenting very elegant 
examples of same in ~\cite{GEA:oddmoduli} and
\cite{GEA:multiRR}.  Since the appearance of those papers, many authors have
presented multisum Rogers-Ramanujan type identities, with a particular emphasis
on infinite families of results, as in Andrews' discovery of the ``Bailey chain"
~(\cite[p. 28 ff]{GEA:qs}, \cite{GEA:multiRR}).  I will not
be presenting such infinite families here as, in light of~\cite{GEA:multiRR}, it
is a routine exercise to imbed any Rogers-Ramanujan type identity  
in such an infinite family.   

There is much current interest in new Bailey pairs and innovations with 
Bailey chains (cf. \cite{GEA:btlcc}, \cite{bw}, \cite{bis}, \cite{jl},
\cite{sow:50}, \cite{sow:ext}).  We shall show in this paper that some
very appealing Rogers-Ramanujan type identities are still to be found 
that are actually derivable from Bailey's original ideas
(combined in some cases with $q$-hypergeometric transformations 
due to Verma and Jain).  For example, 
\begin{gather}
  \sum_{n,r\geqq 0} \frac{q^{2n^2+3r^2+4nr}}{(-q;q)_{n+r} (-q;q)_{n+2r}
    (q;q)_n (q;q)_r}
 = \frac{(q^4,q^5,q^{9};q^{9})_\infty}{(q^2;q^2)_\infty},\label{ex1}\\
\sum_{n,r\geqq 0} \frac{q^{n^2+3r^2+4nr}(-q;-q)_{2n+2r}  }
  {(q^2;q^2)_{2n+2r} (q^2;q^2)_r (q^2;q^2)_n } =
\frac{(q^6,q^8,q^{14};q^{14})_\infty (q^2;q^4)_\infty}{(q;q)_\infty},\\
\sum_{n,r\geqq 0}
  \frac{ q^{3n^2+6nr+6r^2} (q;q)_{3r} }{(q^3;q^3)_{2r} (q^3;q^3)_r (q^3;q^3)_n}
 = \frac{(q^7,q^8,q^{15};q^{15})_\infty}{(q^3;q^3)_\infty},\\
\sum_{n,r\geqq 0} \frac{q^{2n^2 + 3r^2 + 4nr}}
  {(-q;q)_{2n+2r} (q;q^2)_r (q;q)_r (q^2;q^2)_n }
= \frac{(q^{14},q^{16},q^{30};q^{30})_\infty}{(q^2;q^2)_\infty},\\
\sum_{n,r\geqq 0}\frac{q^{n^2+2nr+3r^2} (-q;q^2)_{n+r}}{(q^2;q^2)_n (q^2;q^2)_r (q^
2;q^4)_r}
=\frac{(q^{16},q^{20},q^{36};q^{36})_\infty (q^2;q^4)_\infty}{(q;q)_\infty},\\
\sum_{n,r\geqq 0} \frac{ q^{n^2+2r^2 } }{(q;q^2)_n (q^2;q^2)_r
(q^2,q^4)_r (q;q)_{n-2r} }
= \frac{(q^{28}, q^{32},q^{60};q^{60})_\infty}{(q;q)_\infty}.\label{ex6}
\end{gather}

Furthermore, we note that the double sum identities 
I present here do not arise as merely ``one level up" in the 
standard Bailey chain from well-known
single-sum identities.  

After reviewing the necessary background material in \S\ref{background},
I define the ``standard multiparameter Bailey pair," 
(SMPBP) in \S\ref{smpbp} and
demonstrate that all of the Bailey pairs presented by Bailey in~\cite{wnb1}
and~\cite{wnb2} may be viewed as special cases of the SMPBP.

In \S\ref{sc}, I derive new Bailey pairs as special cases of the SMPBP, and finally
in \S\ref{rrt}, I present a collection of new Rogers-Ramanujan type identities
which are consequences of the Bailey pairs from \S\ref{sc}.   In \S 5, I conclude
by relating this paper's results to Slater's list~\cite{ljs2} and 
suggesting a possible direction of further research.

\subsection{Background}\label{background}
In an effort to understand the mechanism which allowed Rogers to discover the
Rogers-Ramanujan identities and other identities of similar type, 
Bailey discovered
that the underlying engine was quite simple indeed.  This engine was named
the ``Bailey transform" by Slater~\cite[\S 2.3]{ljsb}.

\begin{thm}[The Bailey Transform]
If \[\beta_n = \sum_{r=0}^n \alpha_r u_{n-r} v_{n+r} \] and
   \[ \gamma_n = \sum_{r=n}^\infty \delta_r u_{r-n} v_{r+n},\]
then 
\[\sum_{n=0}^\infty \alpha_n \gamma_n = \sum_{n=0}^\infty \beta_n \delta_n.\]
\end{thm}
  
   Bailey remarks~\cite[p. 1]{wnb2} that ``[t]he proof is almost trivial" and indeed
the proof merely involves reversing the order of summation in a double series.
 
  Curiously, Bailey never uses the Bailey transform in this general form.
He immediately specializes $u_{n} = 1/(q;q)_n$, $v_{n} = 1/(aq;q)_n$, and
\[\delta_n = \frac{(\rho_1;q)_n (\rho_2;q)_n (q^{-N};q)_n q^n}
{(\rho_1 \rho_2 q^{-N} a^{-1};q)_n}.\]  This, in turn, forces
\[ \gamma_n = \frac{(aq/\rho_1;q)_N (aq/\rho_2; q)_N (-1)^n (\rho_1;q)_n (\rho_2;q)_n
(q^{-N};q)_n (aq/\rho_1\rho_2)^n q^{nN-\binomial{n}{2}}}
{(aq;q)_N (aq/\rho_1 \rho_2;q)_N (aq/\rho_1;q)_n (aq/\rho_2;q)_n (aq^{N+1};q)_n} \]

Modern authors normally refer to the $\alpha$ and $\beta$ appearing in the Bailey
transform under the aforementioned specializations of $u$, $v$, and $\delta$ as a
``Bailey pair."
\begin{defn}
A pair of sequences $\left(\alpha_n (a,q),\beta_n(a,q)\right)$ is called a 
\emph{Bailey pair} if
for $n\geqq 0$, 
   \begin{equation} \label{BPdef}
      \beta_n (a,q) = \sum_{r=0}^n \frac{\alpha_r(a,q) }{(q;q)_{n-r} (aq;q)_{n+r}}.
   \end{equation}
\end{defn}
  
In~\cite{wnb1} and~\cite{wnb2}, Bailey proved the fundamental
result now known as ``Bailey's lemma" (see also~\cite[Chapter 3]{GEA:qs}), which
is actually just a consequence of the Bailey transform:
\begin{thm}[Bailey's Lemma]
If $(\alpha_r (a,q), \beta_j (a,q))$ form a Bailey pair, then
\begin{gather} 
  \frac{1}{ (\frac{aq}{\rho_1};q)_n ( \frac{aq}{\rho_2};q)_n} 
  \sum_{j\geqq 0} \frac{ (\rho_1;q)_j (\rho_2;q)_j 
     (\frac{aq}{\rho_1 \rho_2} ;q)_{n-j}}
     {(q;q)_{n-j}} \left( \frac{aq}{\rho_1 \rho_2} \right)^j \beta_j(a,q) \nonumber\\
 = \sum_{r=0}^n \frac{ (\rho_1;q)_r (\rho_2;q)_r}
    { (\frac{aq}{\rho_1};q)_r (\frac{aq}{\rho_2};q)_r (q;q)_{n-r} (aq;q)_{n+r}}
    \left( \frac{aq}{\rho_1 \rho_2} \right)^r \alpha_r(a,q). \label{BL}
\end{gather}
\end{thm}

\section{The Standard Multiparameter Bailey Pair}\label{smpbp}
Since the sequence $\beta_n$ is completely determined for any $\alpha_n$ 
by (\ref{BPdef}), all we need to do is define the standard multiparameter
Bailey pair via the $\alpha_n$ as 
\begin{gather} \label{alphadef} 
    \alpha^{(d,e,k)}_{n} (a,b,q) := \\ 
      \left\{  \begin{array}{ll}
       \displaystyle{\frac{a^{(k-d+1)r/e} q^{(k - d + 1)d r^2/e} 
       (a^{1/e} q^{2d/e};q^{2d/e})_r (a^{1/e};q^{d/e})_r}
       {b^{r/e} {(a^{1/e} b^{-1/e} q^{d/e}; q^{d/e})_r}
        (a^{1/e} ;q^{2d/e})_r (q^{d/e};q^{d/e})_r }, }
            &\mbox{if $n= dr$,} \\
      0,                                      &\mbox{otherwise,}
              \end{array} \right. \nonumber
    \end{gather}
and the corresponding $\beta^{(d,e,k)}_n (a,b,q)$ will be determined
by~\eqref{BPdef}.
Of course, the \emph{form} in which $\beta^{(d,e,k)}_n$ is 
\emph{presented} depends on which $q$-hypergeometric
transformation or summation formula is employed. 
The mathematical interest lies in the fact that elegant Rogers-Ramanujan type identities
 will arise for many choices of $d$, $e$, and $k$, as we shall see in 
\S\ref{rrt}.
    
  \begin{rem} \label{b0inf} An easy calculation reveals that
   \begin{equation}
     \lim_{b\to 0} \alpha^{(d,e,k)}_n (a,b,q) 
     = \lim_{b\to\infty} \alpha^{(d,e,k-1)}_n (a,b,q).
     \end{equation}
 \end{rem}

 \begin{rem}
  In all derivations of Rogers-Ramanujan type identities, 
Bailey lets $b\to 0$ or
  $b\to\infty$.  In light of Remark~\ref{b0inf}, it will be sufficient to let 
  $b\to 0$ from this point forward.  Also, it will be convenient to 
replace $a$ by $a^e$ and $q$ by $q^e$ throughout.  Thus, in practice,
rather than (\ref{alphadef}), we will instead only need to consider the 
somewhat more managable

\begin{gather}  \label{alpha0def}
    \alpha^{(d,e,k)}_{n} (a^e,0,q^e) := \\
        \left\{  \begin{array}{ll}
     \displaystyle{\frac{
                    (-1)^r a^{(k-d)r} q^{(dk - d^2 + \frac d2) r^2 - \frac d2r}
                    (a q^{2d};q^{2d})_r (a;q^d)_r}
        {(a ;q^{2d})_r (q^d;q^d)_r }, }
            &\mbox{if $n= dr$,} \\
      0,                                      &\mbox{otherwise.}
              \end{array} \right. \nonumber
    \end{gather}
\end{rem}
We now proceed to find the corresponding $\beta_n$ for the SMPBP:
 
 \begin{eqnarray*}
  &&\beta^{(d,e,k)}_{n} (a^e,0,q^e)\\
 &=& \sum_{s=0}^n \frac{ \alpha^{(d,e,k)}_s (a^e,0,q^e) }{(q^e;q^e)_{n-s} 
    (a^e q^e ; q^e )_{n+s}} 
\alpha_{d,k,m}(a^e,0,q^e)\\
    &=& \frac{1}{(q^e ; q^e)_{n} (a^e q^e ; q^e)_{n}} 
      \sum_{s=0}^n \frac{(-1)^s  q^{ens - \frac e2 s^2 +\frac e2 s} 
       (q^{-en};q^e)_s}{(a^e q^{e(n+1)} ; q^e)_s} 
      \alpha^{(d,e,k)}_{s} (a^e,0,q^e) \\
      &=& \frac{1}{(q^e ; q^e)_{n} (a^e q^e ; q^e)_{n}} 
      \sum_{r=0}^{\lfloor n/d \rfloor} 
      \frac{(-1)^{dr}  q^{endr - \frac {e d^2}{2} r^2 +\frac {ed}{2} r } 
       (q^{-en};q^e)_{dr} }{(a^e q^{e(n+1)} ; q^e)_{dr} } 
      \alpha^{(d,e,k)}_{dr} (a^e,0,q^e) \\
      &=& \frac{1}{(q^e ; q^e)_{n} (a^e q^e ; q^e)_{n}} 
      \sum_{r=0}^{\lfloor n/d \rfloor} 
      \frac{(a,q^d\sqrt{a}, - q^d\sqrt{a}; q^d)_r 
       (q^{-en};q^e)_{dr}  }
      {(q^d, \sqrt{a},-\sqrt{a}; q^d)_r (a^e q^{e(n+1)} ; q^e)_{dr} }  \\
      & & \qquad\qquad\qquad\qquad\qquad \times
      (-1)^{(d+1)r} a^{(k-d) r}   
      q^{(2k-2d-ed+1) \frac d2 r^2 + (e-1)\frac {d}{2} r + endr }.
\end{eqnarray*}
 Note that 
  \begin{equation} \label{split}
  (q^{-en}; q^e)_{dr} = \prod_{i=0}^{d-1}  (q^{-en + ei};q^{de})_r,
  \end{equation}
 and that each factor in the right hand side of (\ref{split}) can be factored into
 a product of $e$ factors
 \begin{equation*}
 (q^{-en+ei}; q^{de} )_r = \prod_{j=1}^e ( \xi_e^j q^{-n+i}; q^d)_r
\end{equation*}
where $\xi_e$ is a primitive $e$-th root of unity, and 
that the complementary denominator factor 
$(a^e q^{e(n+1)} ; q^e)_{dr}$ can be split and factored analogously into a
product of $ed$ rising $q$-factorials.

   Furthermore, $q^{(2k-2d-ed+1) \frac d2 r^2}$ can be written as a limiting case of
a product of $| 2k-ed-2d+1|$ rising $q$-factorials.
For example, supposing that $k=5$, $e=1$, $d=2$, one can write $q^{5r^2}$
as a limit as $\tau\to 0$ of the product of $|2k-2d-ed+1|  = 5$ 
rising $q$ factorials:
\[ (-1)^r q^{5r^2} = \lim_{\tau\to 0} \tau^{5r} (q/\tau;q^2)_r^5.\]
Thus $\beta^{(d,e,k)}_n (a^e,0,q^e) $ can be seen to be a product of  
${\left( (q^e ; q^e)_n (a^e q^e ; q^e)_n \right)^{-1}}$ and a limiting case of 
a very-well poised ${}_{t+1}\phi_{t}$, where \[t=ed+|2k-ed-2d+1|+2,\] and the
basic hypgeometric series ${}_{p+1}\phi_{p}$ is defined by
\begin{equation*}
  \hgs{p+1}{p}{a_1,a_2,\dots,a_{p+1}}{b_1,b_2,\dots,b_p}{q}{z} =
  \sum_{r=0}^\infty \frac{ (a_1,a_2,\dots,a_{p+1};q)_r}
     {(q,b_1,b_2,\dots,b_p;q)_r} z^r.   
\end{equation*}
It will also be convenient to use the standard abbreviation
\begin{equation*}
\W{p+1}{p}{a; a_3, a_4, \dots, a_{p+1}; q ; z} :=
\hgs{p+1}{p}{a, qa^{\frac 12}, -qa^{\frac 12}, a_3, a_4, \dots, a_{p+1}}
  {a^{\frac 12}, -a^{\frac 12}, \frac{aq}{a_3}, \frac{aq}{a_4}, \dots, 
   \frac{aq}{a_{p+1}}} {q}{z}.
\end{equation*}

\section{Consequences of the SMPBP} \label{sc}
\subsection{Results of Rogers, Bailey, Dyson, and Slater}
 Bailey presented five Bailey pairs in total; he lists them as (i) through (v) on pages 5--6
of~\cite{wnb2}.  These Bailey pairs arise as the following specializations of the
SMPBP:
\begin{itemize}
  \item $(d,e,k)=(1,1,2)$ with $b\to 0$ is equivalent to Bailey's (i),
  \item $(d,e,k)=(1,2,2)$ is equivalent to (ii),
  \item $(d,e,k)=(1,3,2)$ with $b\to 0$ is equivalent to (iii),
  \item $(d,e,k)=(2,1,2)$ is equivalent to (iv), and
  \item $(d,e,k)=(3,1,4)$ with $b\to 0$ is equivalent to (v).
\end{itemize}

Bailey appears to have limited himself to these cases since these 
are the only ones where $\beta^{(d,e,k)}_n (a^e,0,q^e)$ 
is representable as a finite 
product, and thus its insertion into the left hand side of (\ref{BL}) 
will result in a 
\emph{single}-fold sum.  

The point I wish to emphasize here is that since 
$\beta^{(d,e,k)}_n (a^e,0,q^e)$ 
is a finite product times a very well poised basic hypergeometric series,
as long as one is willing to consider multisums,
it is \emph{a priori} plausible that elegant Rogers-Ramanujan type identities
may be derivable for \emph{any} triple $(d,e,k)$ of positive integers, 
as long as one has in hand an appropriate $q$-hypergeometric transformation 
formula.

 Certain specializations of $(d,e,k)$ yield classical Bailey pairs. 
The following table summarizes the
best known classical results which follow from the SMPBP.  
The letter-number codes
in the ``Bailey pair" column refer to the codes used by Slater in her 
two papers
\cite{ljs1} and \cite{ljs2}.  The reason that a single specialization 
of $(d,e,k)$ may 
correspond to more than one of Slater's Bailey pairs is that she 
chose to specialize
$a$ before performing the required $q$-hypergeometric summation or
transformation, and thus listed what corresponds to our $a=1,q,q^2,$ etc., and 
linear combinations thereof, as different Bailey pairs.  
By substituting the Bailey
pairs into various limiting cases of Bailey's lemma, a variety of 
classical identities
result.  The parenthetical numbers in the third column refer to the numbers in
Slater's list~\cite{ljs2}. 

\begin{tabular}{|r|c|l|}
\hline
$(d,e,k)$& Bailey pairs  & RR type identities\\
\hline
(1,1,1)& H17& Euler's pentagonal number theorem (1) \\
(1,1,2)& B1, B3 & Rogers-Ramanujan (18, 14);\\
&  & \qquad\qquad G\"ollnitz-Gordon (36, 34); \\
&  & \qquad\qquad Lebesgue (8, 12)\\
(1,2,2)& G1--G3 & Rogers-Selberg (31--33); Rogers' mod 5 (19, 15)\\
(1,3,2)&      & Bailey's mod 9 identities (41--43)\\
(2,1,2) & C5, C7 & Rogers' mod 10 identities (46, 44)\\
(2,1,3) & C1--C4  & Rogers' mod 14 identities (59--61);\\ 
&  &  \qquad\qquad Rogers' mod 20 (79)\\
(3,1,4) & J1--J6 &  Dyson's mod 27 identities (90--93);\\ 
&  & \qquad\qquad Slater (71--78, 107--116)\\
\hline
\end{tabular}

\begin{rem}
In addition to the classical Rogers-Ramanujan type identities mentioned above,
certain identities discovered more recently can also be derived from the SMPBP.
In particular, the Verma-Jain mod 17 identities~\cite[pp. 247--248, (3.1)--(3.8)]{vj1}
arise from $(d,e,k)=(1,6,3)$, the 
Verma-Jain mod 19 identities~\cite[pp. 248--250, (3.9--3.17)]{vj1} from 
$(1,6,4)$, the Verma-Jain mod 22 
identities~\cite[pp. 250--251, (3.18)--(3.22)]{vj1}
from $(2,1,5)$, an identity of Ole Warnaar related to the modulus 
11~\cite[p. 246, Theorem 1.3; $k=4$]{sow:borwein} from $(2,2,3)$, 
a mod 13 identity of George Andrews~\cite[p. 280, (5.8); $k=1$]{GEA:multiRR}
from $(2,2,4)$,
and finally Dennis Stanton's mod 11 
identity~\cite[p. 65, (6.4)]{ds}
from $(1,2,4)$.
\end{rem}

\subsection{New Bailey pairs}
\begin{rem} In~\cite{avs:rrt}, I presented a number of
Rogers-Ramanujan type identities that arise from what may now be considered
the $b=0$, $e=1$ case of the SMPBP.  The
interested reader is invited to consult~\cite{avs:rrt} for additional examples.
\end{rem}

We now consider specializations of $(d,e,k)$ in 
\begin{equation}\label{SMPBPdef}
\beta^{(d,e,k)}_n (a^e,0,q^e) =  
\sum_{s=0}^n \frac{ \alpha^{(d,e,k)}_s (a^e,0,q^e) }{(q^e;q^e)_{n-s} 
    (a^e q^e ; q^e )_{n+s}},
 \end{equation}
which lead to new identities. 

The most crucial step in each case will be the transformation of the very well poised series
via one of the following known formulas:
\begin{trans}[Watson's $q$-analog of Whipple's Theorem]\cite{gnw},
   \cite[p. 360, equation (III.18)]{gr2} 
 \begin{gather}
   \W{8}{7}{a; b,c,d,e,q^{-n};q,\frac{a^2q^{n+2}}{bcde}}\nonumber\\
    = \frac{(aq,aq/de;q)_n}{(aq/d,aq/e;q)_n}
  \hgs{4}{3}{aq/bc, d, e, q^{-n}}{aq/b,aq/c,deq^{-n}/a}{q}{q}\label{WqW}
  \end{gather}
\end{trans}
is used to establish (\ref{BP123}), (\ref{BP131}), (\ref{BP133}), (\ref{BP142}),
and (\ref{BP143}). 
\begin{trans}[The first Verma-Jain ${}_{10}\phi_9$ transformation]
\cite[p. 232, (1.3)]{vj1}, \cite[p. 97, (3.10.4)]{gr2}
  \begin{gather}
  \W{10}{9}{a; b,x,-x,y,-y,q^{-n},-q^{-n}; q, -\frac{a^3 q^{2n+3}}{bx^2 y^2}}
\nonumber\\
    = \frac{(a^2 q^2, a^2 q^2/x^2 y^2; q^2)_n}{(a^2 q^2/x^2, a^2 q^2/y^2;q^2)_n}
    \hgs{5}{4}{q^{-2n},x^2,y^2,-aq/b, -aq^2/b}
    {x^2 y^2 q^{-2n}/a^2, a^2 q^2/b^2,-aq,-aq^2}{q^2}{q^2}\label{VJ1}
    \end{gather}
\end{trans} 
is used to establish (\ref{BP141}), (\ref{BP143}), (\ref{BP163}), and (\ref{BP164}).

\begin{trans}[The second Verma-Jain ${}_{10}\phi_9$ transformation]
\cite[p. 232, (1.4)]{vj1}
  \begin{gather}
  \W{10}{9}{a; b,x,xq,y,yq,q^{1-n},q^{-n}; q; \frac{a^3 q^{2n+3}}{bx^2 y^2} }
\nonumber\\ 
   = \frac{ (aq, aq/xy; q)_n}{(aq/x, aq/y;q)_n}
    \hgs{5}{4}{x,y,\sqrt{aq/b},-\sqrt{aq/b},q^{-n}}
    {\sqrt{aq},-\sqrt{aq},aq/b,xyq^{-n}/a}{q}{q}\label{VJ2}
    \end{gather}
\end{trans}
is used to establish (\ref{BP215}), (\ref{BP222}),  (\ref{BP225}), and (\ref{BP417}).

\begin{trans}[The first Verma-Jain ${}_{12}\phi_{11}$ transformation]
\cite[p. 232, (1.4)]{vj1}
  \begin{gather}
   \W{12}{11}{a; x,\omega x,\omega^2 x,
       y,\omega y,\omega^2 y, q^{-n},\omega q^{-n}, \omega^2q^{-n}; 
       q; -\frac{a^4 q^{3n+4}}{x^3 y^3} }\nonumber\\
   = \frac{(a^3 q^3, \frac{a^3 q^3}{x^3 y^3}; q^3)_n}
       {(\frac{a^3 q^3}{x^3}, \frac{a^3 q^3}{y^3};q^3)_n}
     \hgs{6}{5}{q^{-3n},x^3,y^3,aq, aq^2,aq^3}
    {(aq)^{\frac 32}, -(aq)^{\frac 32}, a^{\frac 32}q^3,-a^{\frac 32}q^3, 
      \frac{x^3 y^3 q^{-3n}}{a^3}}
    {q^3}{q^3},\label{VJ3}
    \end{gather}
\end{trans} where $\omega = \exp(2\pi i/3)$, 
is used to establish~(\ref{BP135}).

\begin{trans}[The second Verma-Jain ${}_{12}\phi_{11}$ transformation]
\cite[p. 232, (1.5)]{vj1}
  \begin{gather}
    \W{12}{11}{ a; x,xq,xq^2,y,yq,yq^2,q^{2-n},q^{1-n},q^{-n}; q, 
       \frac{a^4 q^{3n+3}}{x^3 y^3} }\nonumber\\
    = \frac{ (aq, aq/xy; q)_n}{(aq/x, aq/y;q)_n}
    \hgs{6}{5}{\sqrt[3]{a}, \omega\sqrt[3]{a}, \omega^2\sqrt[3]{a},x,y,q^{-n}}
    {\sqrt{a},-\sqrt{a},\sqrt{aq},-\sqrt{aq},xyq^{-n}/a}{q}{q}\label{VJ4}
    \end{gather}
\end{trans}
is used to establish~(\ref{BP327}) and (\ref{BP337}).

Finally,
\begin{trans}[Transformation of a very-well-poised ${}_8\phi_7$]
\cite[p. 76, equation (3.4.7), reversed]{gr2}
\begin{gather}
 \W{8}{7}{a; y^{\frac 12},-y^{\frac 12}, (yq)^{\frac 12},-(yq)^{\frac 12},x;
    q; \frac{a^2 q}{y^2 x}}
    = \frac{(aq, \frac{a^2 q}{y^2};q)_\infty}
    { (\frac{aq}{y}, \frac{a^2 q}{y};q)_\infty}
   \hgs{2}{1}{y,\frac{xy}{a}}{\frac{aq}{x}}{q}{\frac{a^2 q}{y^2 x}}. 
  \label{gr8721}
\end{gather}
\end{trans} is used to establish (\ref{BP223}) and (\ref{BP224}).

  Let $(d,e,k) = (1,2,3)$.  Then
 \begin{eqnarray*}
 && \beta^{(1,2,3)}_n (a^2,0,q^2) \\ &=&  
\sum_{r=0}^n \frac{ \alpha^{(1,2,3)}_r (a^2,0,q^2) }{(q^2;q^2)_{n-r} 
    (a^2 q^2 ; q^2 )_{n+r}}\\
 &=& \frac{1}{(q^2 ; q^2)_{n} (a^2 q^2 ; q^2)_{n}} 
      \sum_{r=0}^n \frac{(-1)^r  q^{2nr - r^2 +  r} 
       (q^{-2n};q^2)_r}{(a^2 q^{2(n+1)} ; q^2)_r} 
      \alpha^{(1,2,3)}_{r} (a^2,0,q^2) \\
      &=& \frac{1}{(q^2 ; q^2)_{n} (a^2 q^2 ; q^2)_{n}} 
      \sum_{r=0}^{n} 
      \frac{(a,q\sqrt{a}, - q\sqrt{a}, q^{-n},-q^{-n}; q)_r  a^{2r}   
      q^{2nr + \frac 32 r^2 + \frac r2 } }
      {(q, \sqrt{a},-\sqrt{a}, aq^{n+1},-aq^{n+1}; q)_r }  \\   
 &=& \frac{1}{(q^2 ; q^2)_{n} (a^2 q^2 ; q^2)_{n}} \\
 & & \qquad\times \lim_{\tau\to 0}
  \hgs{8}{7}{a,q\sqrt{a}, - q\sqrt{a}, q^{-n},-q^{-n},q/\tau,1/\tau,q/\tau}
  {\sqrt{a},-\sqrt{a},a q^{n+1},-a q^{n+1},\tau a,\tau a q,\tau a}{q}{-a^2\tau^3
  q^{2n}}\\
  &=& \lim_{\tau\to 0}
  \frac{(aq;q)_n (\tau^2 a;q)_n}{(\tau a;q)_n (\tau aq;q)_n 
  (q^2 ; q^2)_{n} (a^2 q^2 ; q^2)_{n}}\\ 
  &&\qquad\times
     \hgs{4}{3}{-\tau a q^n, q/\tau, 1/\tau, q^{-n}}
       {-aq^{n+1},\tau a, q^{1-n}{a\tau^2}}
  {q}{q} \mbox{(by~(\ref{WqW}))}\\
  &=& \frac{1}{(-q;q)_n} \sum_{r\geqq 0} 
  \frac{a^r q^{r^2}}{(q;q)_r (q;q)_{n-r} (-aq;q)_{n+r}}
 \end{eqnarray*}
 
 Thus, we have established
 \begin{equation}\label{BP123}
\beta^{(1,2,3)}_n (a^2,0,q^2) = \frac{1}{(-q;q)_n} \sum_{r\geqq 0} 
  \frac{a^r q^{r^2}}{(q;q)_r (q;q)_{n-r} (-aq;q)_{n+r}}.
 \end{equation}
 
 Via analogous calculations, one can establish each of the following.  
   
 \begin{equation}\label{BP124}
\beta^{(1,2,4)}_n (a^2,0,q^2) =\sum_{r\geqq 0} 
  \frac{a^{2r} q^{2r^2}}{(-aq;q)_{2r} (q^2;q^2)_{r} (q^2;q^2)_{n-r}}.
 \end{equation} 
 
  \begin{eqnarray} 
 \beta^{(1,3,1)}_n (a^3,0,q^3) &=& 
  \frac{(-1)^n a^{-n} q^{-\binomial{n+1}{2}} (q;q)_n}
 {(q^3;q^3)_n (aq;q)_{2n} }\nonumber\\
&&\qquad\times \sum_{r\geqq 0} 
  \frac{(-1)^r q^{\binomial{r+1}{2} - nr} (aq;q)_{n+r} (aq;q)_{2n+r}}
  {(a^3 q^3;q^3)_{n+r} (q;q)_{r} (q;q)_{n-r}  }. \label{BP131}
 \end{eqnarray}
 
 \begin{equation}\label{BP133}
 \beta^{(1,3,3)}_n (a^3,0,q^3) = \frac{(aq;q)_n }{(aq;q)_{2n} (q^3;q^3)_n}
 \sum_{r\geqq 0} 
  \frac{a^{r} q^{r^2}  (aq;q)_{2n+r} (aq;q)_{n+r}}
  {(a^3 q^3;q^3)_{n+r} (q;q)_{r} (q;q)_{n-r}  }.
 \end{equation}
 
 \begin{equation}\label{BP135}
 \beta^{(1,3,5)}_n (a^3,0,q^3) = 
 \sum_{r\geqq 0} \frac{a^{3r} q^{3r^2} (aq;q)_{3r} }
 {(a^3 q^3;q^3)_{2r} (q^3;q^3)_r (q^3;q^3)_{n-r}}
 \end{equation}
 
 \begin{equation}\label{BP141}
 \beta^{(1,4,1)}_n (a^4,0,q^4) = \frac{(-1)^n q^{2n^2}}{(-a^2 q^2;q^2)_{2n}}
  \sum_{r\geqq 0} \frac{ q^{3r^2-4nr} }{(q^2;q^2)_r (-aq;q)_{2r} (q^4;q^4)_{n-r}}.
 \end{equation} 
 
 \begin{equation}\label{BP142}
 \beta^{(1,4,2)}_n (a^4,0,q^4) = \frac{i^n q^{n^2} (iq;q)_n (q;q)_n}
 {(q^4;q^4)_n (iaq;q)_{2n}^2}
 \sum_{r\geqq 0} 
  \frac{(-i)^r q^{r^2-2nr} (iaq;q)_{2n+r} }
  {(q;q)_r (-iaq;q)_{n+r} (-aq;q)_{n+r} (q;q)_{n-r} (iq;q)_{n-r} }, 
 \end{equation} 
  (where here and throughout, $i=\sqrt{-1}$).
 \begin{equation}\label{BP143}
 \beta^{(1,4,3)}_n (a^4,0,q^4) = \frac{(iq;q)_n (q;q)_n}
 {(q^4;q^4)_n (iaq;q)_{2n}^2} 
 \sum_{r\geqq 0} 
  \frac{a^{r} q^{r^2} (iaq;q)_{2n+r} }
  {(q;q)_r (-iaq;q)_{n+r} (-aq;q)_{n+r} (q;q)_{n-r} (iq;q)_{n-r} }.
 \end{equation} 
 
 \begin{equation}\label{BP144}
 \beta^{(1,4,4)}_n (a^4,0,q^4) = \frac{1}{(-a^2 q^2;q^2)_{2n}}
  \sum_{r\geqq 0} \frac{ a^{2r} q^{2r^2} }{(q^2;q^2)_r (-aq;q)_{2r} (q^4;q^4)_{n-r}}.
 \end{equation} 
 
 \begin{equation}\label{BP163}
 \beta^{(1,6,3)}_n (a^6,0,q^6) = \frac{1}{(a^6 q^6;q^6)_{2n}}
  \sum_{r\geqq 0} \frac{ (-1)^r a^{2r} q^{3r^2} (a^2 q^2;q^2)_{3n-r} }
  {(q^2;q^2)_r (-aq;q)_{2r} (q^6;q^6)_{n-r}}. 
  \end{equation}
 
  \begin{equation}\label{BP164}
 \beta^{(1,6,4)}_n (a^6,0,q^6) = \frac{1}{(a^6 q^6;q^6)_{2n}}
  \sum_{r\geqq 0} \frac{ a^{2r} q^{2r^2} (a^2 q^2;q^2)_{3n-r} }
  {(q^2;q^2)_r (-aq;q)_{2r} (q^6;q^6)_{n-r}} .
  \end{equation} 
 
 \begin{equation}\label{BP215}
 \beta^{(2,1,5)}_n (a,0,q) = \sum_{r\geqq 0} \frac{a^r q^{r^2}}
 {(q;q)_r (aq;q^2)_r (q;q)_{n-r}}. 
 \end{equation}
 
 \begin{equation}\label{BP222}
\beta^{(2,2,2)}_n (a^2,0,q^2) =\frac{(-1)^n q^{n^2}}{(-aq;q)_{2n} }
\sum_{r\geqq 0} 
  \frac{(-1)^r q^{\frac 32 r^2 - \frac 12 r - 2nr}}
{(aq;q^2)_{r} (q;q)_{r} (q^2;q^2)_{n-r}}.
 \end{equation} 
 
 \begin{equation}\label{BP223}
\beta^{(2,2,3)}_n (a^2,0,q^2) =
\frac{(aq^2;q^2)_n}{(a^2 q^2;q^2)_{2n} }\sum_{r\geqq 0} 
  \frac{a^{r} q^{2nr}}{(q^2;q^2)_{r}  (q^2;q^2)_{n-r}}.
 \end{equation} 
 
 \begin{equation}\label{BP224}
\beta^{(2,2,4)}_n (a^2,0,q^2) =
\frac{(aq^2;q^2)_n}{(a^2 q^2;q^2)_{2n} }\sum_{r\geqq 0} 
  \frac{a^{r} q^{2r^2}}{(q^2;q^2)_{r}  (q^2;q^2)_{n-r}}.
 \end{equation}
 
 \begin{equation}\label{BP225}
\beta^{(2,2,5)}_n (a^2,0,q^2) =
\frac{1}{(-a q;q)_{2n} }\sum_{r\geqq 0} 
  \frac{a^{r} q^{r^2}}{(q;q)_{r} (aq;q^2)_r  (q^2;q^2)_{n-r}}.
 \end{equation} 
 
 \begin{equation}\label{BP327}
 \beta^{(3,2,7)}_n (a^2,0,q^2) = 
 \frac{1}{(-aq;q)_{2n}} \sum_{r\geqq 0} \frac{a^r q^{r^2} (a;q^3)_r}
 {(a;q)_{2r} (q;q)_r (q^2;q^2)_{n-r}}.
\end{equation}

 \begin{equation}\label{BP337}
\beta^{(3,3,7)}_n (a^3,0,q^3)=
\frac{1}{(a^3 q^3;q^3)_{2n}} \sum_{r\geqq 0} 
\frac{ a^r q^{3n^2 + r} (a;q^3)_r (aq;q)_{3n-r}}
     { (q;q)_r (a;q)_{2r} (q^3;q^3)_{n-r}}
\end{equation}

\begin{equation}\label{BP417}
\beta^{(4,1,7)}_n (a,0,q) = 
\frac{1}{(aq;q^2)_n} \sum_{r\geqq 0} \frac{a^r q^{2r^2}}
{(q^2;q^2)_r (aq^2;q^4)_r (q;q)_{n-2r}}.
\end{equation}

\section{A list of Rogers-Ramanujan-Bailey type identities}\label{rrt}  
For easy reference, we restate the Bailey Lemma with the SMPBP inserted:
\begin{thm}[Bailey's Lemma]
\begin{gather} 
  \frac{1}{ (\frac{a^e q^e}{\rho_1^e};q^e)_N ( \frac{a^e q^e}{\rho_2^e};q^e)_N} 
 \nonumber \\ 
 \qquad \times \sum_{j\geqq 0} \frac{ (\rho_1^e;q^e)_j (\rho_2^e;q^e)_j 
     (\left( \frac{a q}{\rho_1 \rho_2} \right)^e ;q^e)_{N-j}}
     {(q^e;q^e)_{N-j}} \left( \frac{aq}{\rho_1 \rho_2} \right)^{ej} 
     \beta^{(d,e,k)}_j (a^e, b^e, q^e) \nonumber\\
 = \sum_{r=0}^{\lfloor N/d \rfloor} \frac{ (\rho_1^e;q^e)_{dr} (\rho_2^e;q^e)_{dr}}
    {  ( \left( \frac{aq}{\rho_1} \right)^e;q^e)_{dr} 
        ( \left( \frac{a q}{\rho_2}\right)^e;q^e)_{dr} 
    (q^e;q^e)_{N-dr} (a^e q^e ; q^e)_{N+dr}}
    \left( \frac{a q}{\rho_1 \rho_2} \right)^{der} \nonumber\\ 
  \qquad \times  \alpha^{(d,e,k)}_{dr} (a^e, b^e, q^e). \label{SMPBL}
\end{gather}
\end{thm}

Note that the rather general identities presented by Bailey as
equations (6.1)--(6.4) on page 6 of \cite{wnb2}, from which all of the other
identities in~\cite{wnb1} and \cite{wnb2} can be derived, are simply cases
$(d,e,k) = (1,2,2), (1,3,2), (2,1,2)$, and $(3,1,4)$ respectively,  of (\ref{SMPBL}).

I have compiled a list of Rogers-Ramanujan type identities which I believe to be new.   
Each is a direct consequence of (\ref{SMPBL}), with parameters specialized as indicated 
in brackets.  
Note that the final form of the sum side for many of the identities 
was obtained only after reversing the order of
summation and replacing $n$ by $n+r$.  Of course, in order to obtain each infinite
product representation, Jacobi's triple product 
identity~\cite[p. 15, equation (1.6.1)]{gr2}
is applied to the right hand side after $a$ is specialized.  

This list is by no means exhaustive; I have merely chosen some examples to illustrate
the power of the SMPBP.
\begin{rem} The identity which arises from inserting a given
$(\alpha^{(d,e,k)}_n(a,q),\beta^{(d,e,k}_n(a,q))$ into \eqref{SMPBL} in
the case where $\rho_1,\rho_2, N\to\infty$, and $a=1$ is just one of 
a set of $d(e-1)+k$ identities.  The other identities can be found via
a system of $q$-difference equations.  This phenomenon is explored in
\cite{avs:qdiff}.  
\end{rem} 

\begin{equation} \label{ATNS123}
\sum_{n,r\geqq 0} 
  \frac{ q^{n^2+2nr+2r^2} (-q;q^2)_{n+r}}
  {(-q;q)_{n+r} (-q;q)_{n+2r} (q;q)_n (q;q)_r}
 = \frac{(q^3,q^4,q^7;q^7)_\infty (-q;q^2)_\infty}{(q^2;q^2)_\infty}
\end{equation} [ $N,\rho_1\to\infty$, $\rho_2=-\sqrt{q}$,  $a=1$, $b\to 0$, 
$(d,e,k)=(1,2,3)$.]

\begin{equation} \label{PNS131}
\sum_{n,r\geqq 0} 
  \frac{ (-1)^n q^{\frac 52 n^2- \frac 12 n+4nr+2r^2} 
  (q;q)_{n+r}  (q;q)_{n+2r} (q;q)_{2n+3r}}
  {(q^3;q^3)_{n+r} (q^3;q^3)_{n+2r} (q;q)_n (q;q)_r (q;q)_{2n+2r} }
 = \frac{(q^3,q^4,q^7;q^7)_\infty}{(q^3;q^3)_\infty}
\end{equation} [ $N,\rho_1,\rho_2\to\infty$, $a=1$, $b\to 0$, $(d,e,k)=(1,3,1)$.]

\begin{gather} 
\sum_{n,r\geqq 0} 
  \frac{ (-1)^n q^{2 n^2 -n+ 2nr + r^2} 
  (-q^3;q^6)_{n+r} (q^2;q^2)_{n+r}  (q^2;q^2)_{n+2r} (q^2;q^2)_{2n+3r}}
  {(q^6;q^6)_{n+r} (q^6;q^6)_{n+2r} (q^2;q^2)_n (q^2;q^2)_r (q^2;q^2)_{2n+2r} }
  \nonumber\\  \label{ATNS131}
 = \frac{(q^3,q^5,q^8;q^8)_\infty (-q^3;q^6)_\infty}{(q^6;q^6)_\infty}
\end{gather} [ $N,\rho_1\to\infty$, $\rho_2=-\sqrt{q}$, $a=1$, $b\to 0$, 
$(d,e,k)=(1,3,1)$.]

\begin{equation} \label{PNS123-2}
\sum_{n,r\geqq 0} 
  \frac{ q^{2n^2+2n+4nr+3r^2+3r}}{(-q;q)_{n+r} (-q;q)_{n+2r+1} (q;q)_n (q;q)_r}
 = \frac{(q,q^8,q^9;q^9)_\infty}{(q^2;q^2)_\infty}
\end{equation} [ $N,\rho_1,\rho_2\to\infty$,  $a=q$, $b\to 0$, $(d,e,k)=(1,2,3)$.]

\begin{equation} \label{PNS123}
\sum_{n,r\geqq 0} 
  \frac{ q^{2n^2+4nr+3r^2}}{(-q;q)_{n+r} (-q;q)_{n+2r} (q;q)_n (q;q)_r}
 = \frac{(q^4,q^5,q^9;q^9)_\infty}{(q^2;q^2)_\infty}
\end{equation} [ $N,\rho_1,\rho_2\to\infty$,  $a=1$, $b\to 0$, $(d,e,k)=(1,2,3)$.]

\begin{equation} \label{ATNS124}
\sum_{n,r\geqq 0} 
  \frac{ q^{n^2+2nr+3r^2} (-q;q^2)_{n+r}}{(-q;q)_{2r}  (q^2;q^2)_n (q^2;q^2)_r}
 = \frac{(q^4,q^5,q^{9};q^{9})_\infty (-q;q^2)_\infty}{(q^2;q^2)_\infty}
\end{equation} [ $N,\rho_1\to\infty$, $\rho_2=-\sqrt{q}$, $a=1$, $b\to 0$, 
$(d,e,k)=(1,2,4)$.]

\begin{equation}\label{PNS141-2}
\sum_{n,r\geqq 0}
\frac{(-1)^{n+r} q^{6n^2+8nr+5r^2+4n+4r}}{(-q^2;q^2)_{2n+2r+1} (q^2;q^2)_r 
(-q;q)_{2r+1} (q^4;q^4)_n}
= \frac{(q,q^8,q^9;q^9)_\infty}{(q^4;q^4)_\infty}
\end{equation}
[ $N,\rho_1,\rho_2\to\infty$,  $a=q$, $b\to 0$, $(d,e,k)=(1,4,1)$.]

\begin{equation}\label{PNS141}
\sum_{n,r\geqq 0}
\frac{(-1)^{n+r} q^{6n^2+8nr+5r^2}}{(-q^2;q^2)_{2n+2r} (q^2;q^2)_r (-q;q)_{2r}
(q^4;q^4)_n}
= \frac{(q^4,q^5,q^9;q^9)_\infty}{(q^4;q^4)_\infty}
\end{equation}
[ $N,\rho_1,\rho_2\to\infty$,  $a=1$, $b\to 0$, $(d,e,k)=(1,4,1)$.]

\begin{equation} \label{ATNS222}
\sum_{n,r\geqq 0} 
  \frac{ (-1)^n q^{2n^2+2nr+\frac 32 r^2-\frac r2} (-q;q^2)_{n+r}}
  {(-q;q)_{2n+2r}  (q;q)_r 
  (q;q^2)_r (q^2;q^2)_n}
 = \frac{(q^4,q^{6},q^{10};q^{10})_\infty (-q;q^2)_\infty}{(q^2;q^2)_\infty}
\end{equation} [ $N,\rho_1\to\infty$, $\rho_2=-\sqrt{q}$,
 $a=1$, $b\to 0$, $(d,e,k)=(2,2,2)$.]   
 
 \begin{equation}\label{ATNS141}
\sum_{n,r\geqq 0}
\frac{(-1)^{n+r} q^{8n^2+8nr+6r^2} (-q^4;q^8)_{n+r}}
{(-q^4;q^4)_{2n+2r} (q^4;q^4)_r (-q^2;q^2)_{2r}(q^8;q^8)_n}
= \frac{(q^4,q^6,q^{10};q^{10})_\infty (-q^4;q^8)_\infty}{(q^8;q^8)_\infty}
\end{equation}
[ $N,\rho_1\to\infty$, $\rho_2=-\sqrt{q}$, $a=1$, $b\to 0$, $(d,e,k)=(1,4,1)$.]

\begin{equation} \label{PNS223-2}
\sum_{n,r\geqq 0} 
  \frac{  q^{n^2+2r^2+3nr+2n+3r} (q;q)_{n+r+1} }{(q;q)_{2n+2r+2}  
  (q;q)_r (q;q)_n}
 = \frac{(q,q^{10},q^{11};q^{11})_\infty}{(q;q)_\infty}
\end{equation} [ $N,\rho_1,\rho_2\to\infty$, $a=q$, $b\to 0$, $(d,e,k)=(2,2,3)$],    

\begin{equation} \label{PNS223}    
\sum_{n,r\geqq 0} 
  \frac{  q^{n^2+2r^2+3nr} (q;q)_{n+r} }{(q;q)_{2n+2r}  (q;q)_r (q;q)_n}
 = \frac{(q^5,q^6,q^{11};q^{11})_\infty}{(q;q)_\infty}
\end{equation} [ $N,\rho_1,\rho_2\to\infty$, $a=1$, $b\to 0$, $(d,e,k)=(2,2,3)$]
, 
\emph{(due to S. O. Warnaar~\cite[p. 247, Theorem 1.3 with $k=4$]{sow:borwein}).}

\begin{equation} \label{PNS124-2}
\sum_{n,r\geqq 0} 
  \frac{ q^{2n^2+2n+4nr+4r^2+4r}}{(-q;q)_{2r+1}  (q^2;q^2)_n (q^2;q^2)_r}
 = \frac{(q,q^{10},q^{11};q^{11})_\infty}{(q^2;q^2)_\infty}
\end{equation} [ $N,\rho_1,\rho_2\to\infty$, $a=q$, $b\to 0$, $(d,e,k)=(1,2,4)$.]

\begin{equation} \label{PNS133-2}
\sum_{n,r\geqq 0}
  \frac{ q^{3n^2+6nr+4r^2+3n+4r} (q;q)_{n+r}  (q;q)_{2n+3r+1} (q;q)_{n+2r+1}}
  {(q^3;q^3)_{n+r} (q^3;q^3)_{n+2r+1} (q;q)_n (q;q)_r (q;q)_{2n+2r+1} }
 = \frac{(q,q^{10},q^{11};q^{11})_\infty}{(q^3;q^3)_\infty}
\end{equation} [ $N,\rho_1,\rho_2\to\infty$, $a=q$, $b\to 0$, 
$(d,e,k)=(1,3,3)$.]

\begin{equation} \label{PNS133}
\sum_{n,r\geqq 0} 
  \frac{ q^{3n^2+6nr+4r^2} (q;q)_{n+r}  (q;q)_{2n+3r} (q;q)_{n+2r}}
  {(q^3;q^3)_{n+r} (q^3;q^3)_{n+2r} (q;q)_n (q;q)_r (q;q)_{2n+2r} }
 = \frac{(q^5,q^6,q^{11};q^{11})_\infty}{(q^3;q^3)_\infty}
\end{equation} [ $N,\rho_1,\rho_2\to\infty$, $a=1$, $b\to 0$, $(d,e,k)=(1,3,3)$.]

\begin{gather} 
\sum_{n,r\geqq 0} 
  \frac{ i^n q^{5n^2+8nr+4r^2} (iq;q)_{n+r} (q;q)_{n+r} (iq;q)_{2n+3r} }
  { (q^4;q^4)_{n+r} (iq;q)_{2n+2r}^2 (q;q)_r (-iq;q)_{n+2r} (-q;q)_{n+2r} (q;q)_n
    (iq;q)_n}\nonumber
\\ = \frac{(q^5,q^6,q^{11};q^{11})_\infty}{(q^4;q^4)_\infty} \label{PNS142}
\end{gather} [ $N,\rho_1,\rho_2\to\infty$, $a=1$, $b\to 0$, $(d,e,k)=(1,4,2)$.]

\begin{gather}\label{ATNS144}
\sum_{n,r\geqq 0}
 \frac{q^{2n^2 + 4nr + 4r^2} (-q^2;q^4)_{n+r}}
 {(-q^2;q^2)_{2n+2r} (q^2;q^2)_r (-q;q)_{2r} (q^4;q^4)_n}
= \frac{(q^{5},q^{6},q^{11};q^{11})_\infty (-q^2;q^4)_\infty}
{(q^4;q^4)_\infty}
\end{gather}
 [ $N,\rho_1\to\infty$, $\rho_2=-\sqrt{q}$, $a=1$, $b\to 0$, $(d,e,k)=(1,4,4)$.]

\begin{gather} 
\sum_{n,r\geqq 0}
\frac{(-1)^r q^{3n^2 + 6nr + 6 r^2} (-q^3;q^6)_{n+r} (q^2;q^2)_{3n+2r}}
{(q^6;q^6)_{2n+2r} (q^2;q^2)_r (-q;q)_{2r} (q^6;q^6)_n} \nonumber\\
=\frac{(q^5,q^6,q^{11};q^{11})_\infty (-q^3;q^6)_\infty}{(q^6;q^6)_\infty}
\label{ATNS163}
\end{gather}
[ $N,\rho_1\to\infty$,$\rho_2=-\sqrt{q}$, $a=1$, $b\to 0$, $(d,e,k)=(1,6,3)$.]

\begin{equation} \label{PNS224-2}
\sum_{n,r\geqq 0} 
  \frac{  q^{n^2+2r^2+2nr + 2n + 3r} (q;q)_{n+r+1} }{(q;q)_{2n+2r+2}  (q;q)_r (q;q)_n}
 = \frac{(q,q^{12},q^{13};q^{13})_\infty}{(q;q)_\infty}
\end{equation} [ $N,\rho_1,\rho_2\to\infty$, $a=q$, $b\to 0$, $(d,e,k)=(2,2,4)$.]

\begin{equation} \label{PNS224}
\sum_{n,r\geqq 0} 
  \frac{  q^{n^2+2r^2+2nr} (q;q)_{n+r} }{(q;q)_{2n+2r}  (q;q)_r (q;q)_n}
 = \frac{(q^6,q^7,q^{13};q^{13})_\infty}{(q;q)_\infty}
\end{equation} [ $N,\rho_1,\rho_2\to\infty$, $a=1$, $b\to 0$, $(d,e,k)=(2,2,4)$],
\emph{(due to G. E. Andrews~\cite[p. 280, eq. (5.8) with $k=1$]{GEA:multiRR}).}

\begin{gather} 
\sum_{n,r\geqq 0} 
  \frac{ q^{4n^2+8nr+5r^2} (iq;q)_{n+r} (q;q)_{n+r} (iq;q)_{2n+3r} }
  { (q^4;q^4)_{n+r} (iq;q)_{2n+2r}^2 (q;q)_r (-iq;q)_{n+2r} (-q;q)_{n+2r} (q;q)_n
    (iq;q)_n}\nonumber\\
 = \frac{(q^6,q^7,q^{13};q^{13})_\infty}{(q^4;q^4)_\infty}
\label{PNS143}
\end{gather} [ $N,\rho_1,\rho_2\to\infty$, $a=1$, $b\to 0$, $(d,e,k)=(1,4,3)$.]

\begin{gather}\label{ATNS164}
\sum_{n,r\geqq 0}
 \frac{q^{3n^2 + 6nr + 5r^2} (-q^3;q^6)_{n+r} (q^2;q^2)_{3n+2r}}
 {(q^6;q^6)_{2n+2r} (q^2;q^2)_r (-q;q)_{2r} (q^6;q^6)_n}
= \frac{(q^6,q^7,q^{13};q^{13})_\infty (-q^3;q^6)_\infty}{(q^6;q^6)_\infty}
\end{gather}
 [ $N,\rho_1\to\infty$, $\rho_2=-\sqrt{q}$, $a=1$, $b\to 0$, $(d,e,k)=(1,6,4)$.]

\begin{equation} \label{ATNS223}
\sum_{n,r\geqq 0} 
  \frac{  q^{n^2+3r^2+4nr} (-q;q^2)_{n+r}(q^2;q^2)_{n+r} }
  {(q^2;q^2)_{2n+2r}  (q^2;q^2)_r (q^2;q^2)_n}
 = \frac{(q^6,q^8,q^{14};q^{14})_\infty (-q;q^2)_\infty}{(q^2;q^2)_\infty}
\end{equation} [ $N,\rho_1\to\infty$, $a=1$, $\rho_2=-\sqrt{q}$,
 $b\to 0$, $(d,e,k)=(2,2,3)$.]

\begin{gather} \label{ATNS142}
\sum_{n,r\geqq 0} 
  \frac{ i^n q^{6n^2+8nr+4r^2} (-q^4;q^8)_{n+r} (iq^2;q^2)_{n+r} (q^2;q^2)_{n+r} 
  (iq^2;q^2)_{2n+3r} }
  { (q^8;q^8)_{n+r} (iq^2;q^2)_{2n+2r}^2 (q^2;q^2)_r (-iq^2;q^2)_{n+2r} 
  (-q^2;q^2)_{n+2r} (q^2;q^2)_n
    (iq^2;q^2)_n}\nonumber\\
 = \frac{(q^6,q^8,q^{14};q^{14})_\infty (-q^4;q^8)_\infty}{(q^8;q^8)_\infty}
\end{gather} [ $N,\rho_1,\to\infty$, $\rho_2=-\sqrt{q}$,
$a=1$, $b\to 0$, $(d,e,k)=(1,4,2)$.]

\begin{gather}\label{PNS135-2}
\sum_{n,r\geqq 0} 
  \frac{ q^{3n^2+6nr+6r^2+3n+6r} (q;q)_{3r+1} }
  {(q^3;q^3)_{2r+1} (q^3;q^3)_r (q^3;q^3)_n}
 = \frac{(q,q^{14},q^{15};q^{15})_\infty}{(q^3;q^3)_\infty}
 \end{gather}
  [ $N,\rho_1,\rho_2\to\infty$, $a=q$, $b\to 0$, $(d,e,k)=(1,3,5)$.] 

\begin{gather}\label{PNS135}
\sum_{n,r\geqq 0} 
  \frac{ q^{3n^2+6nr+6r^2} (q;q)_{3r} }{(q^3;q^3)_{2r} (q^3;q^3)_r (q^3;q^3)_n}
 = \frac{(q^7,q^8,q^{15};q^{15})_\infty}{(q^3;q^3)_\infty}
 \end{gather}
  [ $N,\rho_1,\rho_2\to\infty$, $a=1$, $b\to 0$, $(d,e,k)=(1,3,5)$.] 

\begin{gather}\label{PNS144-2}
\sum_{n,r\geqq 0}
 \frac{q^{4n^2 + 8nr + 6r^2+4n+6r}}
 {(-q^2;q^2)_{2n+2r+1} (q^2;q^2)_r (-q;q)_{2r+1} (q^4;q^4)_n}
= \frac{(q,q^{14},q^{15};q^{15})_\infty}{(q^4;q^4)_\infty}
\end{gather}
 [ $N,\rho_1,\rho_2\to\infty$, $a=q$, $b\to 0$, $(d,e,k)=(1,4,4)$.]

\begin{gather}\label{PNS144}
\sum_{n,r\geqq 0}
 \frac{q^{4n^2 + 8nr + 6r^2}}{(-q^2;q^2)_{2n+2r} (q^2;q^2)_r (-q;q)_{2r} (q^4;q^4)_n}
= \frac{(q^7,q^8,q^{15};q^{15})_\infty}{(q^4;q^4)_\infty}
\end{gather}
 [ $N,\rho_1,\rho_2\to\infty$, $a=1$, $b\to 0$, $(d,e,k)=(1,4,4)$.]

\begin{gather} 
\sum_{n,r\geqq 0} 
  \frac{ q^{3n^2+6nr+5r^2} (-q^3;q^6)_{n+r} (q^2;q^2)_{n+r}  
  (q^2;q^2)_{2n+3r} (q^2;q^2)_{n+2r}}
  {(q^6;q^6)_{n+r} (q^6;q^6)_{n+2r} (q^2;q^2)_n (q^2;q^2)_r (q^2;q^2)_{2n+2r} 
  }\nonumber\\ \label{ATNS133}
 = \frac{(q^7,q^9,q^{16};q^{16})_\infty (-q^3;q^6)_\infty}{(q^6;q^6)_\infty}
\end{gather} [ $N,\rho_1\to\infty$, $\rho_2=-\sqrt{q}$, $a=1$, $b\to 0$, 
$(d,e,k)=(1,3,3)$.]

\begin{equation} \label{PNS222-2}
\sum_{n,r\geqq 0} 
  \frac{ (-1)^n q^{3n^2+4n+4nr+\frac 52 r^2+\frac 72 r}}{(-q;q)_{2n+2r+2}  (q;q)_r 
  (q;q^2)_{r+1} (q^2;q^2)_n}
 = \frac{(q^2,q^{16},q^{18};q^{18})_\infty}{(q^2;q^2)_\infty}
\end{equation} [ $N,\rho_1,\rho_2\to\infty$, $a=q^2$, $b\to 0$, $(d,e,k)=(2,2,2)$.]

\begin{equation} \label{PNS222}
\sum_{n,r\geqq 0} 
  \frac{ (-1)^n q^{3n^2+4nr+\frac 52 r^2-\frac r2}}{(-q;q)_{2n+2r}  (q;q)_r 
  (q;q^2)_r (q^2;q^2)_n}
 = \frac{(q^8,q^{10},q^{18};q^{18})_\infty}{(q^2;q^2)_\infty}
\end{equation} [ $N,\rho_1,\rho_2\to\infty$, $a=1$, $b\to 0$, $(d,e,k)=(2,2,2)$.]

\begin{gather} \label{ATNS143}
\sum_{n,r\geqq 0} 
  \frac{ q^{4n^2+8nr+6r^2} (-q^4;q^8)_{n+r} (iq^2;q^2)_{n+r} (q^2;q^2)_{n+r} 
  (iq^2;q^2)_{2n+3r} }
  { (q^8;q^8)_{n+r} (iq^2;q^2)_{2n+2r}^2 (q^2;q^2)_r (-iq^2;q^2)_{n+2r} 
  (-q^2;q^2)_{n+2r} (q^2;q^2)_n
    (iq^2;q^2)_n}\nonumber\\
 = \frac{(q^8,q^{10},q^{18};q^{18})_\infty (-q^4;q^8)_\infty}{(q^8;q^8)_\infty}
\end{gather} [ $N,\rho_1\to\infty$, $\rho_2=-\sqrt{q}$,
$a=1$, $b\to 0$, $(d,e,k)=(1,4,3)$.]
   
\begin{equation} \label{ATNS224}
\sum_{n,r\geqq 0} 
  \frac{  q^{n^2+2nr+3r^2} (-q;q^2)_{n+r} (q^2;q^2)_{n+r} }
  {(q^2;q^2)_{2n+2r}  (q^2;q^2)_r (q^2;q^2)_n}
 = \frac{(q^8,q^{10},q^{18};q^{18})_\infty (-q;q^2)_\infty}{(q^2;q^2)_\infty}
\end{equation} [ $N,\rho_1,\to\infty$, $\rho_2=-\sqrt{q}$,
$a=1$, $b\to 0$, $(d,e,k)=(2,2,4)$.]

\begin{equation}\label{SS215}
\sum_{n,r\geqq 0}
\frac{q^{\frac 12 n^2 + \frac 12 n + nr + \frac 32 r^2+ \frac 12 r} 
(-1;q)_{n+r}}{(q;q)_n (q;q)_r (q;q^2)_r}
=\frac{(q^9,q^9,q^{18};q^{18})_\infty (-q;q)_\infty}{(q;q)_\infty}
\end{equation}
[ $N,\rho_1\to\infty$, $\rho_2=-1$, $a=1$, $b\to 0$, $(d,e,k)=(2,1,5)$.]

\begin{equation} \label{ATNS225}
\sum_{n,r\geqq 0} 
  \frac{  q^{n^2+2r^2+2nr} (-q;q^2)_{n+r} }
  {(-q;q)_{2n+2r}  (q;q^2)_r (q;q)_r (q;q)_n}
 = \frac{(q^{10},q^{12},q^{22};q^{22})_\infty (-q;q^2)_\infty}{(q^2;q^2)_\infty}
\end{equation} [ $N,\rho_1\to\infty$, $\rho_2=-\sqrt{q}$,
 $a=1$, $b\to 0$, $(d,e,k)=(2,2,5)$.]
 
 \begin{gather}\label{ATNS135}
\sum_{n,r\geqq 0} 
  \frac{ q^{3n^2+6nr+9r^2} (-q^3;q^6)_{n+r} (q^2;q^2)_{3r} }
  {(q^6;q^6)_{2r} (q^6;q^6)_r (q^6;q^6)_n}
 = \frac{(q^{11},q^{13},q^{24};q^{24})_\infty (-q^3;q^6)_\infty}
 {(q^6;q^6)_\infty}
 \end{gather}
  [ $N,\rho_1,\to\infty$, $\rho_2=-\sqrt{q}$, $a=1$, $b\to 0$, $(d,e,k)=(1,3,5)$.] 

 \begin{equation} \label{PNS225-2}
\sum_{n,r\geqq 0} 
  \frac{  q^{2n^2+3r^2+4nr+4n+6r} }{(-q;q)_{2n+2r+2}  (q;q^2)_{r+1} (q;q)_r (q^2;q^2)_n}
 = \frac{(q^{2},q^{28},q^{30};q^{30})_\infty}{(q^2;q^2)_\infty}
\end{equation} [ $N,\rho_1,\rho_2\to\infty$, $a=q^2$, $b\to 0$, $(d,e,k)=(2,2,5)$.]
 
 \begin{equation} \label{PNS225}
\sum_{n,r\geqq 0} 
  \frac{  q^{2n^2+3r^2+4nr} }{(-q;q)_{2n+2r}  (q;q^2)_r (q;q)_r (q^2;q^2)_n}
 = \frac{(q^{14},q^{16},q^{30};q^{30})_\infty}{(q^2;q^2)_\infty}
\end{equation} [ $N,\rho_1,\rho_2\to\infty$, $a=1$, $b\to 0$, $(d,e,k)=(2,2,5)$.]

\begin{equation}\label{ATNS215}
\sum_{n,r\geqq 0}
\frac{q^{n^2+2nr+3r^2} (-q;q^2)_{n+r}}{(q^2;q^2)_n (q^2;q^2)_r (q^2;q^4)_r}
=\frac{(q^{16},q^{20},q^{36};q^{36})_\infty (-q;q^2)_\infty}{(q^2;q^2)_\infty}
\end{equation}
[ $N,\rho_1\to\infty$, $\rho_2=-\sqrt{q}$, $a=1$, $b\to 0$, $(d,e,k)=(2,1,5)$.]

\begin{equation}\label{SS417}
\sum_{n,r\geqq 0} \frac{ q^{n(n+1)/2+2r^2 } (-1;q)_n }{(q;q^2)_n (q^2;q^2)_r
(q^2,q^4)_r (q;q)_{n-2r} } 
= \frac{(q^{22}, q^{22},q^{44};q^{44})_\infty (-q;q)_\infty}{(q;q)_\infty}
\end{equation}
[ $N,\rho_1\to\infty$, $\rho_2=-1$, $a=1$, $b\to 0$, $(d,e,k)=(4,1,7)$.]

\begin{equation}\label{PNS417-2}
\sum_{n,r\geqq 0} \frac{ q^{n^2+2r^2 +4n+4r} }{(q;q^2)_{n+2} (q^2;q^2)_r
(q^2,q^4)_{r+1} (q;q)_{n-2r} } 
= \frac{(q^4, q^{56},q^{60};q^{60})_\infty}{(q;q)_\infty}
\end{equation}
[ $N,\rho_1,\rho_2\to\infty$, $a=q^4$, $b\to 0$, $(d,e,k)=(4,1,7)$.]

\begin{equation}\label{PNS417}
\sum_{n,r\geqq 0} \frac{ q^{n^2+2r^2 } }{(q;q^2)_n (q^2;q^2)_r
(q^2,q^4)_r (q;q)_{n-2r} } 
= \frac{(q^{28}, q^{32},q^{60};q^{60})_\infty}{(q;q)_\infty}
\end{equation}
[ $N,\rho_1,\rho_2\to\infty$, $a=1$, $b\to 0$, $(d,e,k)=(4,1,7)$.]

\begin{equation}\label{PNS327}
\sum_{n,r\geqq 0} \frac{q^{2n^2 + 4nr + 3r^2 + 6n + 9r} (q^3;q^3)_r}
{(-q;q)_{2n+2r+3} (q;q)_{2r+2} (q;q)_r (q^2;q^2)_n}
= \frac{(q^3,q^{60},q^{63};q^{63})_\infty}{(q^2;q^2)_\infty}
\end{equation}
[ $N,\rho_1,\rho_2\to\infty$, $a=q^3$, $b\to 0$, $(d,e,k)=(3,2,7)$.]

\begin{equation}\label{PNS337}
\sum_{n,r\geqq 0} \frac{ q^{3n^2 + 4r^2 + 6nr + 9n + 12r}
(q^3;q^3)_r (q;q)_{3n+2r+3}}
{ (q^3;q^3)_{2n+2r+3} (q;q)_r (q;q)_{2r+2} (q^3;q^3)_n }
= \frac{(q^3,q^{78},q^{81};q^{81})_\infty}{(q^3;q^3)_\infty}
\end{equation}
[ $N,\rho_1,\rho_2\to\infty$, $a=q^3$, $b\to 0$, $(d,e,k)=(3,3,7)$.]

\begin{equation}\label{ATNS417}
\sum_{n,r\geqq 0} \frac{ q^{n^2+4r^2 } (-q;q^2)_n}{(q^2;q^4)_n (q^4;q^4)_r
(q^4;q^8)_r (q^2;q^2)_{n-2r} } 
= \frac{(q^{40}, q^{48},q^{88};q^{88})_\infty (-q;q^2)_\infty}{(q^2;q^2)_\infty}
\end{equation}
[ $N,\rho_1\to\infty$, $\rho_2=-\sqrt{q}$, $a=1$, $b\to 0$, $(d,e,k)=(4,1,7)$.]

\section{Conclusion}
As remarked earlier, after discovering his lemma and transform, Bailey considered a total of 
five Bailey pairs, and together with Dyson, derived quite a few identities from them
(\cite{wnb1},\cite{wnb2}).
Soon after Bailey completed his work on Rogers-Ramanujan type identities,
Slater, by her own count, found 94 Bailey pairs, leading
to 130 identities (although both of these totals are admittedly somewhat inflated as 
redundancies exist in both).

 For about a quarter century following the work of Bailey and Slater, a handful of
mathematicians did work related to the Rogers-Ramanujan identities, notably
Henry Alder, George Andrews, David Bressoud, Leonard Carlitz, and Basil Gordon.  
Then around 1980, an explosion of interest in Rogers-Ramanujan occured in the mathematics
and physics communities as connections were found with Lie algebras (thanks to
Jim Lepowsky, Steve Milne, and Robert Wilson) and 
statistical mechanics (through the efforts of 
Rodney Baxter, Alex Berkovich, Barry McCoy, Anne Schilling,
Ole Warnaar and others).  A bit later, as the computer
revolution in mathematics began, important contributions related to Rogers-Ramanujan 
were made by Peter Paule, Axel Riese, Herb Wilf, Doron Zeilberger and others.
Some of these practitioners have improved and extended Bailey's lemma.  For example,
Andrews~\cite{GEA:multiRR} showed how the Bailey lemma is self-replicating,
leading to the so-called ``Bailey chain."  Andrews, Schilling, and Warnaar~\cite{asw:A2}
found an ``$A_2$ Bailey lemma."  Andrews and Berkovich
(\cite{GEA:btlcc}, \cite{ab:WP}) extend the Bailey chain to a ``Bailey tree." 
Further innovations and extensions of the Bailey chain are given by
Berkovich and Warnaar~\cite{bw}; Bressoud, Ismail, and Stanton~\cite{bis};
Jeremy Lovejoy~\cite{jl}; and Warnaar~\cite{sow:ext}. 
Ismail and Stanton obtained Rogers-Ramanujan type identities via tribasic
integration~\cite{is}.  
In~\cite{vps}, V. Spiridonov found an elliptic analog of 
the Bailey chain.  And the list goes on and on.  For additional 
references, see the end notes of Chapter 2 in the new edition of
Gasper and Rahman~\cite{gr2}.

The inspiration for this current paper comes from a desire to 
``return to the basics" and
to gain an understanding of Bailey's contributions via a unification of his 
work.  
Accordingly, I refer to (\ref{SMPBPdef}) as the \emph{standard} 
multiparameter Bailey pair because so many of the classical 
Rogers-Ramanujan type identities are direct consequences
of it.  As noted earlier, all identities in Bailey's papers~\cite{wnb1} 
and \cite{wnb2}
may be derived from the SMPBP.  A rough count indicates that at least 40 
percent of
Slater's list~\cite{ljs2} may be derivable from the SMPBP.  However, it seems
plausible that other multiparameter Bailey pairs, analogous to the SMPBP 
in some way,
could be defined, perhaps accounting for the rest of (or at least large 
portions of the
rest of) 
Slater's list, and incidentally revealing many new 
identities of the Rogers-Ramanujan type.  
This clearly warrants further investigation.

\section*{Acknowledgements}
I thank Christian Krattenthaler for writing the HYPQ Mathematica 
package and making
it freely available.  The use of the HYPQ package was invaluable for 
this project.  Also, I thank Bruce Berndt, Ole Warnaar, and
the anonymous referee for helpful comments, and George Andrews
for his continued encouragement.

\end{document}